\documentclass[11pt]{amsart}

\usepackage[margin=1in]{geometry}
\usepackage{amsmath,amssymb,amsthm,mathtools}
\usepackage{hyperref}
\usepackage[shortlabels]{enumitem}
\usepackage{soul}

\usepackage{xcolor}
\usepackage{todonotes}
\usepackage[all,color]{xy}

\theoremstyle{remark}

\DeclareMathOperator{\lat}{lat}

\title{Banach lattices and phase retrieval: A case study for the use of AI in mathematics}

\author[J. de Dios Pont]{Jaume de Dios Pont}
\address{Center for Data Science, New York University, New York, New York 10011, USA}
\email{jdedios@nyu.edu}

\author[L. Liehr]{Lukas Liehr}
\address{Department of Mathematics, Bar-Ilan University, Ramat-Gan 5290002, Israel}
\email{lukas.liehr@biu.ac.il}

\author[D. Muñoz-Lahoz]{David Muñoz-Lahoz}
\address{Instituto de Ciencias
Matemáticas. Universidad Autónoma de Madrid, Campus de Cantoblanco UAM
28049 Madrid, Spain}
\email{david.munnozl@uam.es}

\author[M. A. Taylor]{Mitchell A. Taylor}
\address{ Department of Mathematics\\
ETH Z\"urich, Ramistrasse 101, 8092 Z\"urich, Switzerland}
\email{mitchell.taylor@math.ethz.ch}

\author[P.~Tradacete]{Pedro Tradacete}
\address{Instituto de Ciencias Matem\'aticas (CSIC-UAM-UC3M-UCM)\\
Consejo Superior de Investigaciones Cient\'ificas\\
C/ Nicol\'as Cabrera, 13--15, Campus de Cantoblanco UAM\\
28049 Madrid, Spain.}
\email{pedro.tradacete@icmat.es}

\begin{document}

\begin{abstract}
The ability of large language models to assist professional mathematicians has been progressing rapidly. Earlier this year, a group of researchers in Banach lattice theory and phase retrieval began incorporating this technology into their research workflows. Facing challenges about the reliability of these models, they also decided to couple the discovery process with Lean verification. Here, we present a case study of how this has led to a more united community and a deeper understanding of our field.
\end{abstract}

\maketitle

\section{Introduction}
The purpose of this note is to give an account of some of the progress made in the last five months in Banach lattice theory and phase retrieval, with a focus on how we have begun to incorporate technology into our research. Our objective will be to explain how Lean, large language models (LLMs), and a dedicated team of mathematicians working in unison can achieve an improved understanding of their field, a more united community, and even uncover new opportunities to connect their research with other areas of STEM. 
\vspace{1em} 

In the following, we briefly explain the three essential ingredients that were combined to achieve a broader understanding of our field of mathematics. Interestingly, this enhanced understanding occurred at many different levels, beginning with the core theory and the overall geometry of our field (as we needed to build a Lean library of useful facts, which required us to revisit the classical theory and better explore how known results depended on each other) and culminating at cutting edge research. 

\subsection{Lean}
Lean \cite{lean4,mathlib2020} is an interactive theorem prover and functional programming language designed for writing mathematically precise definitions, statements, and formal proofs. At its core, Lean is based on dependent type theory, where propositions are represented as types and proofs are terms inhabiting those types. The Lean kernel is responsible for checking that each term has the claimed type. Under the propositions as types correspondence, this means that the Lean kernel checks, using only the foundations of dependent type theory, that a given proof is correct.
\vspace{1em}

Thus, if we have a Lean verified proof of a theorem, and we are sure that the statement of the theorem in Lean code faithfully represents the statement in the way we traditionally write mathematics, then we can be certain that the statement holds in the given axioms of mathematics. 
However, proof checking has traditionally been limited by scale.
Lean proofs must contain all of the details---there is no room for interpretation. Even though Lean does a great job at providing algorithms (called tactics) that can help you fill in terms automatically, writing complex mathematics in Lean is challenging and time-consuming. This is where the coding capabilities of generative AI can be exploited to write Lean at a scale useful for research-level mathematics; this process is described in more detail in Section~\ref{Sect:LLM}.

\subsection{Large language models} 
Large language models are AI systems trained on vast collections of text and code to generate language. Given a prompt, an LLM produces an answer which is a likely continuation in terms of the patterns that it has learned during training. For mathematical research, LLMs can be very helpful for identifying possible connections between fields and producing sketches of proofs. However, there is no mathematical certainty in their arguments. Moreover, although some progress has been made in solving certain problems using LLMs, the vast majority of open research questions cannot be given to the models we used with any expectation that they will solve them correctly.

\subsection{Our community of mathematicians} 

We are a group of researchers in Functional Analysis and Applied Harmonic Analysis who have been incorporating Lean and LLMs into our daily mathematical lives. This note describes some of our experiences over a 4-5 month span. During this period, we used LLMs and Lean as research tools for both mathematical exploration and verification. Throughout the process, we continued our usual mathematical duties, such as selecting meaningful problems, contributing numerous novel ideas, and connecting our work to relevant literature. However,  we also aimed to recognize useful contributions from LLMs, which required digesting and evaluating their outputs and deciding whether it could lead to a deeper understanding of our field. The purpose of this note is to provide a concrete and transparent account of what worked for us, what failed, and how these tools affected  our research and our community. However, we do not wish to claim that any of our practices are optimal or should be incorporated by the wider mathematical community.
\vspace{1em}

Several related articles help to put our work in context. Tao in \cite{tao2025machine} surveys a broad range of machine assisted methods and presents LLMs as tools that can complement expert mathematical judgment. Buzzard, Commelin, and Massot demonstrated in  \cite{buzzard2020formalising}  that researchers in mathematics with no direct background in computer science can become effective in the usage of formal verification methods in a reasonably short amount of time. Davies et al.~\cite{davies2021advancing} describe a  methodology in which machine learning techniques identify patterns that professional mathematicians interpret and convert into conjectures and proofs. The present note brings these themes together with LLM-assisted exploration, Lean verification, subject-specific library building and the careful outreach and training of a mathematical community. 
\vspace{1em}

It is important to emphasize that we are not the only mathematicians incorporating LLMs into their research. By now, there are hundreds of articles on ArXiv that disclose the usage of AI tools. There are also many initiatives such as the DARPA expMATH program and numerous AI4Math startups including Harmonic, Axiom Math and Math, Inc. Various perspectives on AI in mathematics have also been shared at the Future of Mathematics Symposium at Stanford and in white papers such as \cite{avigad2026mathematicians,klowden2026mathematical}. Within these reactions stands the Leiden Declaration on Artificial Intelligence and Mathematics, published on June 2, 2026, and endorsed by the International Mathematical Union. This declaration has called for responsible governance of AI in mathematical research. It is a community statement, which up to now has been signed by around 3000 people \url{https://leidendeclaration.ai/}. It argues that although AI tools may affect many aspects of our job, these tools should not alter core values such as correctness, transparency, honesty, attribution,  ethics, open science, independent verification, human understanding, community values, and research autonomy. In this article, we will describe some of the benefits that we have received by following these protocols, and how they have led not only to an improved understanding of our subject, but also to deeper scientific relationships with our colleagues.  
\vspace{1em}

\subsection{Outline of the article}
In Section~\ref{Section2}, we give a brief introduction to our mathematical field. This section is technical and can be skipped on the first reading. In Sections~\ref {Sect:Formalization} and \ref{Sect:Hotspots}, we recall our experiences using Lean and LLMs in the years 2024 and 2025. Then, in Section~\ref{Section: First auto}, we explain how we recently combined these technologies. A more detailed summary of the progress made in the last five months in phase retrieval and related areas is given in Section~\ref{Sect:Papers+Lib}.
\vspace{1em}

In Section~\ref{Section:BLlib}, we briefly explain our new Banach lattice library, referring the reader to \cite{BanLib} for more information. A discussion on the importance of Mathlib and the shortcomings of our library is presented in Section~\ref{Sect:Mathlib}. An explanation of why we have chosen to formalize our results is given in Section~\ref{Section:Why?}. Finally, Section~\ref{SectionRA} records some of the recent LLM-assisted progress in Banach lattice theory (mainly by students). We remark that this latter section contains several technical terms; the reader is encouraged to ignore all undefined terms and focus on the main message. 
\vspace{1em}

The article concludes with Sections~\ref{Sect:community} and \ref{Sect:Future}, which discuss some of the community aspects of this endeavor, as well as a few future plans.



\section{Banach lattices and phase retrieval}\label{Section2}
In this section, we give some basic background and motivation for our field. The reader who is solely interested in the discussions on Lean and LLMs can skip to Section~\ref{Sect:LLM}. 
\vspace{1em}

The \emph{phase retrieval problem} arises in numerous physical scenarios, ranging from crystallography and speech recognition to quantum mechanics. Mathematically, one begins with an injective linear operator $T:H \to X$ which maps a signal space $H$ to a function space $X$. The operator $T$ usually describes some physical phenomena, where, in practice, one is only able to measure the modulus $|Tx|$ of the signal $x\in H$. A classical example is the \emph{Pauli problem}, where $T: L^2(\mathbb{R}^d)\to L^2(\mathbb{R}^d)\times L^2(\mathbb{R}^d)$ is the mapping $T\psi=(\psi,\mathcal{F}\psi)$. In this example, $|T\psi|$ grants access to $|\psi|^2$ and $|\mathcal{F}\psi|^2,$ which are interpreted as the probability densities of position and momentum of a quantum particle. Pauli's problem asks to what extent such data allows one to recover the original wavefunction $\psi$.
\vspace{1em}

Note that in the phase retrieval problem, there is always an inevitable ambiguity; namely, for any unimodular scalar $\lambda$, we have $|Tx|=|T(\lambda x)|$. Therefore, the best one can hope to do is recover the original signal $x$ up to global phase. This ambiguity is often harmless in practice, as global phase does not usually carry physical meaning. With this in mind, one can frame the uniqueness problem for phase retrieval as follows: Given an injective linear operator $T$, is it true that $|Tx|=|Ty|$ if and only if $x=\lambda y$ for some unimodular scalar $\lambda$?
\vspace{1em}

Other common examples of operators $T$ which are relevant for phase retrieval include the Fourier transform $\mathcal{F}$ and the short-time Fourier transform $V_h$ with window function $h\in L^2(\mathbb{R}^d)$. These phase retrieval problems occur in imaging applications. Indeed, in classical crystallography, if one calls the image one wants to see $f$, then the detector records the power spectrum $|\mathcal{F}f|^2$. Thus, one is tasked with the problem of recovering the original object $f$ from the distorted image $|\mathcal{F}f|^2$. Since the Fourier transform is a unitary operator, this task is clearly impossible, unless one knows a priori information about the object $f$. However, one often does,\footnote{A natural constraint may, for example, be the support properties of $f$.} meaning one may assume that $f\in F\subseteq L^2(\mathbb{R}^d)$ for a proper subspace $F\subseteq L^2(\mathbb{R}^d)$. In this case, the model operator is $T:=\mathcal{F|}_F$.
\vspace{1em}

Since obtaining a priori information on $f$ is often difficult, practitioners have designed experiments that incorporate additional redundancy. For example, in ptychography, a mask or window $h$ is inserted, so that one measures $|\mathcal{F}(fh)|^2$. One then translates the window $h$ (or sometimes the object $f$), resulting in measurements of the form $|\mathcal{F}(fT_xh)|^2$, where $T_x$ denotes the translation operator by $x$. In other words, one now measures $|V_hf|^2$ and wishes to recover $f$. Since $V_h: L^2(\mathbb{R}^d)\to L^2(\mathbb{R}^{2d})$, this experimental setup affords a whole variable of redundancy, which one hopes to use to offset the loss of phase information without needing to impose additional constraints on $f\in L^2(\mathbb{R}^d)$. Mathematically, it immediately connects the phase retrieval problem to the field of \emph{time-frequency analysis}.
\vspace{1em}

The discussion above has been entirely algebraic. However, in practice, there will always be measurement errors, so one does not only desire uniqueness in the phase retrieval problem but also \emph{stability}. In other words, one desires an inequality of the form
\begin{equation}\label{T-stable}
    \inf_{|\lambda|=1}\|x-\lambda y\|\leq C \||Tx|-|Ty|\|, \hspace{5mm} x,y\in H.
\end{equation}
In practice, the operator $T$ is bounded and, moreover, if such an inequality is to hold, then $T^{-1}$ must also be bounded. Therefore, we have the proportionality 
  \begin{equation*}
      \inf_{|\lambda|=1}\|x-\lambda y\|\sim \inf_{|\lambda|=1}\|Tx-\lambda Ty\|.
  \end{equation*}
  Up to changing the constant $C$, we may therefore rewrite \eqref{T-stable} as
\begin{equation}\label{stable}
    \inf_{|\lambda|=1}\|f-\lambda g\|\leq C \||f|-|g|\|, \hspace{5mm} f,g\in E:=T(H).
\end{equation}
The utility of the reformulation \eqref{stable} is that it removes the operator $T$ from the nonlinear inequality and instead uses the operator $T$ to select the subspace $E:= T(H)\subseteq X$ where the nonlinear inequality is desired to hold. It also immediately connects stable phase retrieval to the theory of \emph{Banach lattices}, which, in essence, can be defined as the spaces where all of the symbols in \eqref{stable} make sense.
\vspace{1em}

Banach lattices are the natural abstraction of function spaces such as $L_p(\mu)$, $C(K)$, Orlicz spaces and Lorentz spaces. They are defined as Banach spaces $(X,\leq)$ equipped with a lattice order (so that $|f|:= f\vee (-f):=\sup\{f,-f\}$ exists) satisfying the natural compatibility conditions:
\begin{enumerate}
    \item $f,g,h\in X$, $f\leq g$ implies $f+h\leq g+h$.
    \item $f\in X$, $\lambda\in \mathbb{R}_+$, $f\ge 0$ implies $\lambda f\geq 0$.
    \item $|f|\leq |g|$ implies $\|f\|\leq \|g\|.$
\end{enumerate}
These objects have been intensively studied for decades, but it was only recently noticed that the theory has major implications on the
mathematics of phase retrieval.
\vspace{1em}

Banach lattices are not only modeled after spaces of functions, but
can be truly represented as such. For instance,
Kakutani's $C(K)$-representation theorem guarantees that every Banach
lattice is (in a possibly inequivalent norm) ``locally'' a space of continuous functions on some compact
Hausdorff space $K$. Banach lattices also admit global
representations as closed sublattices of $\ell_\infty $-sums of
$L_1$-spaces. However, in some sense, they are much more than just spaces of
functions. For example, they have a rich categorical structure, as shown by the recent
development of free objects in Banach lattices and related categories (see \cite{OTTT} and the references therein).
\vspace{1em}

A relevant part of the literature on Banach lattices is concerned with the interaction with Banach space theory. However, there are also natural overlaps with vector lattice and order theory, which have a more algebraic flavor. 
The theory of Banach lattices also has a longstanding interaction with operator theory and a number of applications in mathematical economics. The \emph{Positivity} conference series serves as a consolidated meeting point for the community.

\section{Our experience using Lean and Large Language Models}\label{Sect:LLM}

Here, we will give a broad overview of our experiences using LLMs and Lean to improve our mathematical research.
\subsection{Our first steps into formalization}\label{Sect:Formalization}
The first result in Banach lattice theory to be formalized was Kakutani's $C(K)$-representation theorem. This theorem states that any Banach lattice whose unit ball contains a largest element must, in fact, be a space of continuous functions on a compact Hausdorff space $K$. Moreover, in this case, the largest element can be identified with the constant function $\mathbf{1}$. 
\vspace{1em}

Given any positive element $e$ in a Banach lattice $X$, we may define the principal ideal
\begin{equation*}
    I_e:=\{x\in X : \exists \, \lambda\in \mathbb{R} \,\text{with}\, |x|\leq \lambda e\}.
\end{equation*}
This ideal has a natural norm, $\|\cdot\|_e$, defined for $x\in I_e$ as the smallest real number $\lambda\geq 0$ for which $|x|\leq \lambda e$. It is routine to show that this norm is complete. Hence, by Kakutani's theorem, $(I_e,\|\cdot\|_e)$ is lattice isometric to a space of continuous functions on a compact Hausdorff space $K_e$.
\vspace{1em}

The fact that Banach lattices are ``locally" $C(K)$-spaces (albeit in the $\|\cdot\|_e$ norm rather than the original $\|\cdot\|_X$ norm) is fundamental to the theory. In particular, it implies that any inequality between real numbers that involves only finitely many linear and lattice operations will remain true in any Banach lattice. Indeed, if $x_1,\dots, x_n\in X$, then we may define $e:=|x_1|+\cdots+|x_n|$ and view $x_1,\dots ,x_n \in I_e =C(K_e)$. In $C(K_e)$, the order is pointwise, so if an inequality involving only linear and lattice operations is true for real numbers, it is also true in $C(K_e)$, and hence in $X$. 
\vspace{1em}

Kakutani's $C(K)$-representation theorem was formalized by DML on December 19, 2025. The formalization took around one month, did not involve any use of AI and can be accessed via the following link: 
\begin{center}
    https://github.com/davidmunozlahoz/veclat.
\end{center}  
The motivation for DML to learn Lean and formalize this theorem was mainly curiosity. Intrinsic interest was also the main catalyst for JD to learn Lean, who then taught a seminar and a full course on formalization at ETH Zürich in 2025. These lectures were attended by S.~Bertolini and MT. As will be discussed below, formalization is not easy, and having experienced people who knew how to formalize mathematics by hand was essential to our project.

\subsection{Our first experience with LLMs}\label{Sect:Hotspots}
The first experience where an LLM had a positive impact on our research occurred in late 2024, with the paper \cite{pont2025sharp}. Here, our goal was to precisely quantify the failure of the hot spots conjecture by computing the hot spots ratio $S_d$. 
\vspace{1em}

Interesting results from~\cite{mariano2023improved,steinerberger2023upper} show that $S_d\leq 60$ and $\limsup_{d\to\infty}S_d\leq \sqrt{e}$. Motivated by an informal suggestion in \cite{jerison2000hot}, we used GPT to help perform a relatively routine computation using a family of spherical Neumann sieves to give a lower bound on the hot spots ratio. To our surprise, this showed that $\liminf_{d\to\infty}S_d\geq \sqrt{e}$.
\vspace{1em}

In view of the above, it was natural to conjecture that the spherical sieve construction would give the exact value of $S_d$ in every dimension $d$.  After several weeks of thinking about this problem, we asked GPT for some ideas. In a multi-turn conversation, GPT suggested combining three standard ingredients that we thought were too naive to be sharp. Although the execution of the LLM was completely incorrect, we did eventually produce a corresponding upper bound on $S_d$ using exactly these tools.
\vspace{1em}

We emphasize that the above experience was not at all reproducible and we were unable to extract another novel mathematical insight from an LLM until March 2026.
\subsection{Combining LLMs and Lean}\label{Section: First auto}
 The first paper that we autoformalized was \cite{SPRRV}. The main bulk of the formalization was completed in early-mid March, 2026. The motivation for formalizing this result was mainly to see if it was \emph{possible}. However, we learned through trial and error that not only was autoformalization possible, but it can also be highly \emph{instructive}. Indeed, our original proof in \cite{SPRRV} relied on a concentration compactness scheme, which we found to be extremely difficult to formalize. We therefore spent a significant amount of time working with LLMs to make the proof simpler and more quantitative. By the end of the project, we had not only a Lean verification of our theorem, but an alternative, more elementary (though perhaps less illuminating) proof. Both of these proofs were presented in \cite{SPRRV}, as they contained interesting and complementary ideas.
\vspace{1em}

The biggest shift in our perspectives on the usage of LLMs in mathematics came during a visit of JD to ETH Zürich during the week of March~19-25, 2026. During this week, we were working on the $L^2$-stability problem for STFT phase retrieval. This is a problem that we had thought about for a long time in collaboration with several other mathematicians. We had made significant progress on it, but had not yet reached a final solution and had no indication that we would be able to. We therefore tried to see if we could get some new inspiration by discussing the problem with GPT 5.4 Plus. During a multi-turn conversation, GPT offered a largely incorrect execution of a possible strategy to prove the main lemma to which we had reduced the local stability problem at the constant function $\textbf{1}$. However, as mathematicians familiar with the problem, we quickly realized that the overall strategy was sensible and would allow us to complete our proof.  Moreover, using our previous notes, our autoformalization experience from  \cite{SPRRV} and guiding GPT to correctly fill in the missing details, we were able to obtain a Lean verification of the solution in less than 72 hours. Interestingly, a large percentage of the formalization was completed while we were visiting PSI, which contains a leading ptychography lab in Switzerland. As we will discuss in Section~\ref{Sect:Future}, we hope that LLMs will help bridge the gap between physics and mathematics, leading to numerous new and exciting collaborations.
\vspace{1em}

As mathematicians, our job is not only to solve problems, but to find beautiful proofs and to ask natural followup questions \cite{thurston2006proof}. In the first experience mentioned above \cite{SPRRV}, we were able to find a new and instructive proof of our result. In the second experience \cite{bertolini20262}, it was  clear to us that the proof would generalize significantly. Indeed, although our original goal was to prove local stability in the Fock space (identified with the image of the STFT with Gaussian window) at the constant function in one dimension, it was clear that the proof architecture would generalize to prove local stability for all Hermite windows at all elements in the finite span of the canonical basis and in all dimensions. However, it was also clear that the proof would be largely technical and not instructive. We therefore experimented, breaking the generalization into multiple steps, giving the LLMs the Lean-verified proof from the previous step, working with them to get a proof of the next step, Lean verifying that step, and then repeating the strategy at the next level of generality. This process took roughly six weeks, which we estimate to be about the same amount of time that it would have taken us to write the proof carefully by hand. However, as mathematicians with significant experience writing technical proofs, this alternative strategy was both more instructive and more enjoyable for us.

\subsection{Next steps: Formalizing papers and building libraries}\label{Sect:Papers+Lib}
After the above two experiences, we contacted several members of the phase retrieval and Banach lattice communities to see if they would be interested in an experiment to systematically incorporate LLMs into our mathematical lives. In certain cases, this was met with extreme enthusiasm and in other cases with skepticism or concern for the future of mathematics. The large variance in reactions is, of course, completely natural, and we therefore tried to work with people on as personal a level as possible. This ranged from some members using LLMs on a daily basis, others wishing to have a weekly meeting to discuss team progress, and others who preferred to wait until we had strong evidence that this would be a viable new way of doing mathematics and a systematic way to teach them how to properly use it. In all of the above cases, respect for people's opinions was our top priority. With that being said, we will comment below on some of the outcomes of the last few months.
\vspace{1em}

As mentioned above, our first experience with autoformalization was \cite{SPRRV}. Here, we already had an informal proof of our main result, and we wanted to formalize it in Lean. The key lesson was that formalization often forces one to rethink and simplify the proof and in this case it lead to new insights on the problem. We remark that these insights cannot be observed from the theorem statement itself, which did not change. However, mathematics is not just about solving problems. It is also a quest to find deeper understanding, which makes having multiple proofs with different insights fundamentally important. 
\vspace{1em}

Our second experience with autoformalization was \cite{bertolini20262}. This involved a major collaboration with LLMs both at the formalization stage and at the proof discovery stage. This problem was introduced to MT by J.~Ramos, who was able to prove the result for squares. The proof of this is instructive and will likely be submitted for publication soon. However, we could not adapt the method to prove the result in general. Instead, we developed other strategies, important reductions, and began discussing with other people as well. As mentioned above, the key breakthrough occurred when GPT 5.4 offered a new perspective on the problem, incorrectly executed, but observed to be important by us. We then properly incorporated this new insight into the scheme of proof that we had been developing before. Another interesting application of the LLMs was to help generalize the result  -- based on the same ideas -- to a more natural statement. This roughly doubled the length of the proof, but, more importantly, allowed the authors to focus the exposition on the ideas rather than on  the routine, but technical, details.
\vspace{1em}

Another interesting experience occurred in the article \cite{Lukas1}. In this paper, we investigated a problem that we believed to be within reach, but with a solution that lied outside of the expertise of anyone who had worked on the problem previously. Indeed, this was a problem in time-frequency analysis where the solution was naturally obtained by a topological argument involving Chern classes. Our interest in this problem came from our desire to better understand time-frequency analysis and its connection to phase retrieval. Notably, we are not experts in topology, but we felt that this problem had a topological flavor. We therefore worked with GPT 5.4 to obtain our desired counterexample. However, not being complete experts in all of the material, we also wanted to know with absolute certainty that the proof was correct. We therefore autoformalized the result in Lean, which gave us complete confidence that our argument was sound. Moreover, as in \cite{bertolini20262}, we naturally extended the counterexample to its rightful level of generality. This experience highlighted that solutions to certain problems that a priori seem out of reach may sometimes be obtained by using techniques from other areas of mathematics. This leads both to an incredible opportunity to connect fields and learn new mathematics, and to a concern that the proposed execution is indeed correct. In our case, Lean helped ensure that our argument was indeed airtight.
\vspace{1em}

A second issue in \cite{Lukas1} was that formalizing the \emph{statement} of our theorem in Lean in a manner digestible to a non-expert was quite difficult. In fact, in our original theorem, the statement would have likely been impossible for anyone except a Lean expert to understand. Therefore, we spent a lot of time and effort reformulating and refining our theorem. In the end, we found a much simpler way to state the result. Interestingly, the task of properly formalizing and refining the theorem in Lean also led us to discover several new examples, equivalences and structural phenomena related to the underlying theory that we would not have found if we had not put in this effort.
These additional findings will be reported on in a separate article, which will likely involve more difficult Lean statements and definitions. In particular, in this case, we expect that it will be more efficient to carefully formalize the concepts that we are repeatedly using and incorporate them into a reusable library.
\vspace{1em}

The next paper that we highlight in detail is \cite{Lukas2}. This paper was correctly Lean verified -- with almost no guidance -- by LL. This is quite remarkable, as he had no experience with Lean three months ago. Of course, those more experienced with Lean carefully checked that the verification is correct, which was an extremely important part of the process. This success was made possible by a broader effort to help our community to get started with Lean. Indeed, we have been designing tutorials to conveniently set up and teach Lean to our community \cite{leantutorial}.  We also emphasize that translating statements from natural language into Lean is not at all trivial and that the models we used struggled significantly with formalizing definitions. Therefore, a careful, hands-on approach to teaching Lean was an imperative part of our community project. 
\vspace{1em}

As mentioned above, autoformalization \emph{can} lead to improved understanding. However, we wish to acknowledge that  the opposite can also be true. Indeed, we had one rather surprising situation where GPT 5.4 in Codex quickly constructed an ambiguity function with a bounded spectral gap, which we knew had some implications in STFT phase retrieval. The construction was technical and, being skeptical that it was correct, we decided to run it through Lean. Surprisingly, within 48 hours we had a \emph{certificate} of the correctness of the theorem. However, in the spirit of Thurston's insightful article \cite{thurston2006proof}, we did not feel that we had a \emph{proof} of this theorem. To avoid a Faustian deal, it was therefore very important for us to revisit the argument, understand it, make it more general, make it more enlightening, and to track the entirely omitted references to where the ideas came from.
As Klowden and Tao emphasize, formal correctness of an argument captures only the logical validity but it does not communicate to a reader why the proof works, what ideas it is based on, or how the argument might generalize \cite[Section 4.4]{klowden2026mathematical}.
In this vein, we would also like to stress the importance of traditional methods of mathematical communication, such as discussing with colleagues and preparing slides for talks. Indeed, the intent to communicate mathematics often leads to a deeper understanding of the mathematics itself, and in our case JD even discovered a new mechanism behind one of our recent proofs while preparing a talk on it.
\vspace{1em}

We have had several other experiences in addition to the above. Personal and group dynamics will be discussed further in Section~\ref{Sect:community} and experiences (mainly by students) with the use of LLMs as research assistants can be found in Section~\ref{SectionRA}. For now, we emphasize that there have been many more failures than successes. Indeed, although we have focused on the successes above, the majority of the time the LLMs did not know how to make any progress on the problem and just outputted nonsense. Moreover, being limited by context length, the models we used rarely produced anything reliable that was more than 5-10 pages. This made human intuition and guidance imperative, as answering a difficult problem often took months of human-AI discussions. However, as we remarked above, there were a few cases in which the LLMs were able to solve a problem remarkably quickly. Although these proofs never required truly new ideas, it was sometimes striking how well the LLMs were able to answer \emph{certain} problems. As mentioned above, in this case, our job was to understand the proof, put it in its proper generality, do a literature search to find out where the method of proof came from, and then ask and investigate the next natural followup question. 
\subsection{The Banach lattice library}\label{Section:BLlib}
When formalizing papers, our main objective was to get a Lean verification. However, when using LLMs directly, the quality, reusability, and organization of the code can be very poor. The only parts of the code that were intensively reviewed by us were the \emph{statement} of the theorem and the \emph{definitions} it relied on. These were decoupled from the Lean proof, and carefully explained in both the paper and the GitHub repository. We were also careful to pick problems where the statement, together with all definitions, could be formalized in an easily digestible manner. We remark that neither of these approaches scale: One instead wants to build detailed, reusable libraries, so that one can easily verify complicated theorems which rely on many definitions in an efficient manner.
\vspace{1em}

Thus, with this approach in mind, the
    library \cite{BanLib} will serve two main purposes. On the one hand, it will be
    a trusted base of knowledge, in which the community reaches a consensus that the formalized definitions correspond to the
    actual definitions used in the field. As noted above, when
    working on an autoformalization project, it is not uncommon for
    LLMs to make slight changes to the definitions that greatly
    simplify the proofs of the theorems. Thus, careful checking of the
    definitions becomes paramount. This means not only checking the
    definitions that appear in the theorem, but also the definitions that
    these definitions depend upon, and so on. The only way one can
    perform this recursive task efficiently is if, as close to the
    main theorem as possible, definitions only depend on a
    trusted library, and therefore one can be certain that they faithfully
    represent the concepts that they claim to represent.  In particular, we stress that formalization
    merely shifts the burden of trust from the proofs to the statements and definitions. With the inherent scale and
    complexity of mathematics, it is therefore necessary for numerous people to carefully review and certify the material in the library so that it can be trusted.
    \vspace{1em}
     
    On the other hand, libraries serve as code infrastructure to make
    most formalizations not only feasible, but also efficient. If,
    eventually, we want formalization of research to become a
    mathematical standard (at least in our community),
    sustainability and affordability of the autoformalization process
    are critical. No one knows what the cost of LLM tokens will be in the future, but it is in our best interest to keep the ratio of mathematical
    output to token usage as high as possible. One clear way to
    accomplish this is by having shared and reusable code in the form
    of a library.  
    \vspace{1em}

 This is the reason we started the Banach lattice library \cite{BanLib}. This library was coded making an intensive use of LLMs, but with careful human supervision; we could say that it was ``semiautoformalized.'' Although the resulting code cannot be compared to that of purely human-written libraries, the supervision process made sure that the definitions are correct, the approaches to formalizing the different concepts are reasonable, and that the organization of the library is clean, understandable and practical for researchers. 
\vspace{1em}

   At the same time, the sacrifice in code quality had a benefit in the form of scale, to the point where we could use it to naturally formalize three research-level results and constructions. We expect that, as other researchers join the project, they will greatly expand the library in their directions of interest and be able to formalize their own research.
Indeed, we have reasons to believe that this level of scale is possible. The library project was launched at the beginning of April, 2026, and so far DML has been its sole contributor. After two months of part-time work, and starting from the definition of a Banach lattice, we got the library to a point where the formalization of the three research-level results mentioned above, including the solution to an open question that was solved very recently, was possible. Needless to say, our library builds on top of Mathlib; this would not have been possible without the enormous human effort that has been put into writing, reviewing and curating Mathlib, as discussed in the next section.
\vspace{1em}

    As a desired by-product of building the library, we gained new mathematical insights into the theory of Banach lattices, which are reported in detail in \cite[Section 3]{BanLib}. This emphasizes the idea that the process of building the library is not only about producing formally-verified facts, but also about improving our understanding of mathematics.
\subsection{The importance of Mathlib}\label{Sect:Mathlib}
As mentioned above, none of our formalization projects would have been possible without Mathlib. Mathlib is the core reason for our trust in our formalizations,
and this trust stems from the fact that Mathlib is thoughtfully designed, from definitions, to theorems, to proofs. Moreover, the extreme care of the Mathlib  community in formalizing mathematical objects has ensured that their meaning faithfully represents what the mathematical community expects. This, in turn, has ensured that the material in Mathlib may be intuitively and reliably extended to formalize more advanced concepts. Indeed, if such care were not taken, formalizations of advanced topics would become exceedingly messy, and hence more prone to misinterpretation.
\vspace{1em}

Our library is not yet at the level of Mathlib.  Its purpose is to facilitate mathematical research, and for this reason we have had to sacrifice code quality. However, it is an evolving project that we hope
will both expand and improve over time. In this regard, we would like to emphasize the formalization scales suggested in \cite{tao2026}. The formalization of our papers (so far) is at the lowest level of formalization, as the goal is solely to verify correctness. However, we would argue that our library lies somewhere between levels 2
(Publication-level formalization) and 3 (Prototype-level
formalization). The increase in quality of the library compared to that of the papers is to ensure reusability and trust. However, our library leaves much to be desired. In particular, we wish not only for it to expand over time, but also for the core material to be revisited. Indeed, this has already happened, as after proving a new result we realized that we needed to frame certain concepts in the library in the language of locally solid topologies. Importantly, given the infrastructure that we had already built, this was not difficult, as many of the proofs were identical. It was also an instructive experience, as the authors learned about the utility of locally solid topologies from redoing parts of the library at this level of generality.
\vspace{1em}

Complementing the scales of formalization in \cite{tao2026}, we also remark that there are scales of \emph{understanding}. Three of the authors are very familiar with executing technical details in Banach lattice theory, but, through the process of semiautoformalization, we learned a significant amount
about the large scale structure of our field. We believe that the scales of formalization and understanding can be symbiotic. Indeed, the process of carefully formalizing a proof leads to an extreme understanding of the fine details of a small subfield, as it causes one to deeply reflect on the mathematical structures at play. However, semiautoformalizing a library also leads to fundamentally new insights about the overall architecture of one's field. 
\vspace{1em}

We remark that we have also limited our scale compared, e.g., to \cite{rammal2026formalizing}. Our goal was to build a reusable, trusted and affordable research-level library in our specific research field, so we did not pursue utmost scale or speed. Instead, we picked a pace where human understanding and trust in the library could be easily obtained. In the near future, we hope to encourage researchers in our field to extend the library at the same quality level to cover their topics of interest, and concurrently revisit the core of the library to improve code quality, put results in proper generality, and encourage students and researchers to learn more about the fine details of the field. We also emphasize that we are continually learning, and refer the reader to \cite{ilin2026sorries} for additional perspectives on formalization.


\subsection{Why we are formalizing our papers}\label{Section:Why?}
The main reason that we have decided to formalize our papers is to ensure correctness against LLM hallucinations.  Indeed, there have been several instances (see Section~\ref{SectionRA}) where we were able to generate a  proof that looked entirely correct, but with a subtle error -- even in just a single line -- invalidating the whole argument. In many cases, finding the mistake was not trivial: It required a lot of patience, mathematical intuition and skepticism, and, without this, the mistake could have easily been missed.\footnote{Especially since we wanted to believe that the result was true.} Moreover, we have found that  models are getting better at hiding their lies. Indeed, we have had a recent experience where five independent agents certified essentially the same argument, which turned out to be completely wrong.
\vspace{1em}

We also note another experience, where an LLM offered a proof of a major conjecture in Banach lattice theory. In this case, the reasoning of the LLM was \emph{entirely correct}. However, it relied on an external reference, which we later read and found to have a gap. Again, this is an example where skepticism and mathematical intuition allowed us to find the error in the reasoning. However, it also brought us new awareness on potential difficulties in carefully validating LLM-assisted proofs.
\vspace{1em}

Formalizing mathematics in Lean is a \emph{partial} solution to this issue. However, Lean is \emph{not} a black box that outputs the correctness of your theorem. Indeed, formalizing definitions and translations of natural language statements into Lean code requires tremendous skill and precision (in \cite{jiang2026} this is called the \emph{specification gap}: Lean can formally verify the correctness of a statement but it cannot guarantee that the statement faithfully expresses the intended mathematical theorem). Learning how to do this takes time, and has not been easy for our group. To partially remedy this, we have designed tutorials to help people in our community learn Lean, and we have experts available to answer any questions and to certify that Lean translations are correct. Notably, this strategy does not immediately scale, as we can only reliably teach these high standards to a small group of devoted mathematicians.
\vspace{1em}

We have also found that formalizing theorems in Lean gets harder as one gets to more advanced material. Indeed, introductory textbooks usually contain all of the details, and a natural path to formalize the material would be to follow the textbook. However, research papers are far from self-contained and usually rely on dozens of external inputs, so this blueprint completely breaks down. Moreover, a nice way to double check that a theorem is correctly formalized is to make sure that it is \emph{useful}, i.e., when you invoke it as an ingredient in a later proof, it actually serves its purpose. However, research papers often contain results that are not repeatedly used, so this method of double checking breaks down. For this and many other reasons (e.g., incorrect references as mentioned above) partial formalizations can also be misleading.
\vspace{1em}

The above paragraph once again highlights why we have decided to build trusted libraries, as discussed in Section~\ref{Section:BLlib}. We believe that reasonable code quality and excellent organization will be necessary for our community to efficiently verify deep theories and, more importantly, to be able to quickly check that the verification is correct. Building such libraries will require large scale communication, and we believe that this is an aspect where our community approach will shine.
\vspace{1em}

Finally, we mention that the desire to formally verify an important and technically delicate theorem need not be specific to LLM-influenced mathematics. When proposing the Liquid Tensor Experiment \cite{scholze2022liquid}, Scholze explained his motivation as wanting to be sure that one of his most important theorems to date was undoubtedly correct.

\subsection{Our experiences using LLMs as a research assistant}\label{SectionRA}
In the following, we highlight experiences collected from feedback from students\footnote{This information was collected before the release of the most recent models, which have improved in many aspects.} of PT on several ongoing research projects in Banach lattice theory, vector lattices, free constructions, and the homological properties of Banach lattices. They show a consistent pattern: LLMs were often useful as mathematical assistants in the early and intermediate stages of research to propose examples, suggest proofs, identify plausible characterizations, and in some cases produce complete arguments. However, their output required systematic verification. Some of the most useful contributions came with incomplete or incorrect proofs, and several failures were caused by the model relying on false general principles that sounded plausible in the surrounding mathematical context. 
\vspace{1em}


 We begin with a work in progress on free dual Banach lattices. For Banach spaces, the bidual behaves naturally as a free object in the category of dual Banach spaces with adjoint (equivalently weak$^\ast$-to-weak$^\ast$ continuous) operators. In the setting of Banach lattices, the corresponding question is subtler, since lattice homomorphisms are not stable under taking adjoints. Indeed, the adjoint of a lattice homomorphism is generally an interval preserving operator rather than a lattice homomorphism. This led to the problem of determining the Banach lattices which admit a free object in the category of dual Banach lattices with adjoint lattice homomorphisms, and in particular of characterizing when the bidual of the free Banach lattice has this universal property. In this case, the LLM produced a complete and understandable solution after only a small number of iterations. The result is now part of an ongoing project of A.~Avil\'es, E.~Garc\'ia-S\'anchez, and G.~Mart\'inez-Cervantes, building on earlier work on free dual spaces and free Banach lattices \cite{GT}.
\vspace{1em}

A related success arose while studying weak$^\ast$ closed sublattices generated by sets. If $A$ is a subset of the dual of a Banach lattice, one may consider the operation
\[
    A \longmapsto \overline{\lat(A)}^{w^\ast}.
\]
Since lattice operations need not be weak$^\ast$ continuous, the weak$^\ast$ closure of a sublattice need not itself be a sublattice. Previously, examples had been constructed for scattered compact spaces $K$ and subsets $A\subseteq C(K)^\ast$ such that
\[
    Z_1=\overline{\lat(A)}^{w^\ast}
\]
is not a sublattice, while
\[
    Z_2=\overline{\lat(Z_1)}^{w^\ast}
\]
is already a weak$^\ast$ closed sublattice. The open question was how many iterations of this combined operation might be necessary before stabilization. The LLM first produced an example, based on the same underlying ideas, requiring three iterations rather than two. It then suggested a mechanism that appears capable of producing examples requiring an arbitrary ordinal number of steps before stabilization. The details of this construction are still being developed by the authors, but the model's contribution was valuable because it pointed to a broader transfinite phenomenon rather than merely solving a single finite case.
\vspace{1em}

Another substantial group of contributions concerned free products and push-outs of vector lattices, Archimedean vector lattices, and Banach lattices. In an ongoing project by G. Mart\'inez-Fern\'andez and PT, the LLM helped formulate and prove several results about preservation of injectivity in categorical constructions. For example, in a push-out square of vector lattices
\[
\xymatrix{X_1 \ar[r]^{S_1} & PO \\
X_0 \ar[u]^{T_1} \ar[r]_{T_2} & X_2 \ar[u]_{S_2}}
\]
the model suggested a proof that if $T_2$ is injective, then $S_1$ is injective. The proof strategy replaces the Hahn--Banach type argument from \cite{AT} with extension and contraction of prime ideals, followed by a reduction to totally ordered vector spaces.
\vspace{1em}

The LLM also helped with a key injectivity statement for Banach lattice free products. If $T_1:C(K_1)\to X_1$ and $T_2:C(K_2)\to X_2$ are injective lattice homomorphisms, then the induced lattice homomorphism
\[
    T_1\overline{\ast}T_2:C(K_1)\ast C(K_2)\longrightarrow X_1\ast X_2
\]
is injective. The suggested argument used the representation $C(K_1)\ast C(K_2)\simeq C(K_1\ast K_2)$ developed in \cite{MFT} but supplied the missing ingredient together with further structural consequences, including the injectivity of the canonical map from the Archimedean vector lattice free product of a family of Banach lattices into their Banach lattice free product, and the order density of $X_1\ast_{\mathrm{AVL}}X_2$ inside $X_1\ast X_2$.
\vspace{1em}

In the same circle of questions, the LLM was useful for generating examples and counterexamples. It helped identify that the vector lattice free product of two Archimedean vector lattices need not be Archimedean, even though free vector lattices themselves are Archimedean. The obstruction was translated into a question about open polyhedral cone covers of closed cones in $\mathbb{R}^n$; the cases where such covers fail yield non-Archimedean free products. These examples illustrate one of the strongest uses of LLMs in this work: They are good at proposing concrete test cases, pushing definitions through specific situations, and suggesting counterexamples when a statement is too optimistic. This was particularly helpful in all of the cases mentioned above, as the questions being pursued were rather ``new" and the authors were still in the process of gaining intuition for them. In particular, we emphasize that the models that we employed in this case study were much less useful in shedding light on longstanding problems or on problems that we had been stuck on for a long time.
\vspace{1em}

A further line of work involved homological properties of Banach lattices, especially extension and splitting phenomena. In Banach space theory, triviality of exact sequences is closely related to extension properties. In the category of Banach lattices, however, one must distinguish between an extension property and a splitting property. For $\lambda>1$, a Banach lattice $E$ has the $\lambda$-extension property if for every lattice isometric interval preserving embedding $i:Y\to X$ and every lattice homomorphism $T:Y\to E$, there is a lattice homomorphism $S:X\to E$ such that $S i=T$ and $\|S\|\leq \lambda\|T\|$. On the other hand, $E$ has the splitting property if every lattice isometric interval preserving embedding $i:E\to X$ has image which is a projection band in $X$. The interval preserving assumption ensures that the image is an ideal, allowing quotients to remain within the relevant category.
\vspace{1em}

The LLM helped clarify that the extension property implies the splitting property, but not conversely. It produced a correct proof that $C[0,1]$ does not have the extension property, while the researchers already knew that $C[0,1]$ has the splitting property. More precisely, it generated an example of a lattice homomorphism into $C[0,1]$ that cannot be extended. This argument was later generalized by the authors to spaces $C(K)$ where $K$ contains non-trivial convergent sequences.
\vspace{1em}

The model also suggested useful characterizations. For the splitting property, it proposed the criterion that a Banach lattice has the splitting property if and only if every truncation is internal. The statement turned out to be correct, although the proof supplied by the LLM had substantial gaps. The final proof was written independently and was simpler than the model's version. For the extension property, the researchers asked whether monotone completeness implied the extension property and suspected that the converse might fail. The model instead answered that no counterexample should exist, because the two properties are equivalent. This statement was correct, although again the proof supplied by the model was not. The researchers eventually proved the equivalence using some of the ideas suggested by the model. These results form part of an advanced work in progress \cite{BSTT}.
\vspace{1em}

In the above case study, one recurring issue was that the LLM would tend to use false facts that are close to true statements in the surrounding theory. For instance, in the work involving weak$^\ast$ convergence, it repeatedly asserted that every weak$^\ast$ convergent net in the dual of a Banach space is norm bounded, invoking the uniform boundedness principle. This is clearly false, as there are weak$^\ast$ convergent nets in dual spaces that are unbounded in norm. The mistake is particularly dangerous because the analogous statement for sequences, or for pointwise bounded families under the hypotheses of uniform boundedness, can make the argument superficially plausible.
\vspace{1em}

Another repeated failure concerned AM-spaces (spaces which satisfy $\||f|\vee|g|\|=\|f\|\vee\|g\|)$. The model repeatedly used the false assertion that every AM-space is isomorphic, as a Banach space or as a lattice, to a space $C(K)$. This occurred in the study of the splitting property for AM-spaces.  Other recurrent mistakes followed the same pattern. In the study of vector lattices, the model often tried to use lattice homomorphisms into $\mathbb{R}$ in situations where such homomorphisms need not exist. A further issue we find important to mention was incomplete citation: On many occasions, the arguments proposed by the model can be identified with known facts, but the proper source and authorship were entirely omitted. In addition, the model sometimes invoked results from the literature while omitting essential hypotheses. These errors are beyond merely bibliographic, as omitting a hypothesis can invalidate an entire proof.
\vspace{1em}

Overall, these experiences suggest that LLMs in this context are most reliable when used as engines for mathematical exploration rather than as final authorities. They can quickly generate proof strategies, test examples, identify alternative formulations, and suggest plausible characterizations. They can also sustain long searches that would be tedious to perform manually. At the same time, every AI-generated argument needs to be checked line by line, or better yet in Lean. The failures show that the large language models used above may repeatedly import false general facts, blur distinctions between nearby categories, and cite theorems without their full hypotheses. Used critically, however, they functioned as productive research assistants that accelerated the generation of ideas and uncovered results that have now become part of ongoing research papers.


\section{Our community}\label{Sect:community}
As mentioned above, it was important to us that this was a community project. Indeed, the past few months have seen a tremendous increase in the number of quality interactions between members of our group. It was also important for us to give many talks and to build tutorials that could potentially benefit other interested members of the Positivity, Functional Analysis and Applied Harmonic Analysis communities. The interactions over the past few months were facilitated by numerous group chats, weekly online meetings and in-person interactions. In particular, MT is grateful to J.~Serra who funded his travel during this period and allowed for vital in-person interactions with various members of the group.
\vspace{1em}

A primary goal of this project was to benefit students. At the beginning of the project, we insisted that students complete any projects that did not involve AI, so that they could receive full credit for their excellent contributions. Moreover, in the papers where they did use AI, we asked them to carefully disclose the process, identify what is \emph{their} contribution, verbally explain the arguments to us, and maintain an extremely high standard of quality. Otherwise, we let them explore on their own. Although the sample size is small, we found that students adapted more quickly to the new technology than  senior members and had some of the best attitudes of all members of the group. This can be clearly seen by the progress shown in Section~\ref{SectionRA}, which hints that properly mentored students can learn faster, strive for higher standards, and achieve deeper understanding than ever before.
\vspace{1em}

Large scale communication and organization was also a very important aspect of our group's success. In particular, it was vital to respect people's ideas, give people space to attack their problems, and to keep morale high. It was also important for students to investigate problems that they cared about (rather than just collecting random problems an LLM happened to know how to solve) and to publish results that they are proud of. Although organizing such a large group during a time of change was not easy, in general it brought the community closer.
\vspace{1em}

Finally, we remark on the adoption of Lean into our mathematical workflows. For the most part, the community has shown a lot of enthusiasm for formalizing their results, although some members have not yet been convinced that this is important. The main roadblock is the high entry point\footnote{Though this has been decreasing with time.}. For this reason, we have spent a lot of time teaching, and also tried to present the Lean statements in our papers in a way that was easier for mathematicians to digest rather than optimized for coding.\footnote{We did not do this in the library, of course.}  In our small sample size, we have noticed that those using Lean experienced an increase in morale and creativity. Indeed, reading pages of LLM text can be exhausting, and it is much nicer to share a result with colleagues that you know is correct. Similarly, exploring out-of-the-box ideas becomes easier when you have multiple ways to later verify correctness.
\vspace{1em}

We also note that we have a wonderful group of Bachelor's students at ETH Zürich who took JD's class on Lean and have been developing new insights on AI-assisted formalization under the guidance of S.~Bertolini and J.~Serra. Their enthusiasm has been palpable, and we have been extremely lucky to benefit from their numerous insights.




\section{Future endeavors}\label{Sect:Future}
Since it is impossible to predict the future, here we will focus on the concrete plans that we have already committed to. Our first objective is to increase communication. We have already been asked to give numerous talks, workshops, and minicourses on AI for math. In addition, we will organize our own conference on AI for analysis at IMPA in early 2027, and DML will give a Lean minicourse at our conference on phase retrieval in September 2026 \url{https://www.mathc.rwth-aachen.de/phaseretrieval/}. JD will also give a course on formalization at NYU in Fall 2026 which will bring together graduate students from Mathematics, Computer Science and Data Science, and MT has proposed dissertations on formalization in analysis for up to four students at Oxford University.
\vspace{1em}

In addition to the above, we are also in the process of improving organization.
 This includes not only organization at the personal level but also in designing new infrastructure to organize relevant results and partial progress. As we have seen, one strength that LLMs have is that they can invoke facts from all over mathematics. In general, the literature is poorly organized, as we have traditionally rewarded advancing the forefront of knowledge much more than we have collecting, digesting, and properly presenting known mathematics. However, with AI tools, better organization is possible and canonizing mathematical knowledge into usable databases and Lean libraries seems both doable at scale and fundamentally important. 
\vspace{1em}

Otherwise, we wish to further develop our field and expand our scientific knowledge. As mentioned above, it is not uncommon for LLMs to solve a problem by invoking a fact that you did not know. This then forces you to expand your mathematical horizons and learn new fields. In the course of using LLMs, we have also needed to take on many jobs which are atypical for mathematicians. Indeed, many of the difficulties we have overcome can be better classified as software engineering problems. Moreover, in our research, we have for the first time used tools such as interval arithmetic to certify a proof, and we expect that using large scale numerics will become common for testing hypotheses. Thus, we expect our mathematical lives to blend somewhat with computer science and engineering and have a more direct impact on the other sciences. More concretely, we have never had more confidence that we may be able to help answer real world problems in phase retrieval and related areas, such as the dynamic range problem. In particular, we hope that the discussions with physicists and engineers at our conference in September will stimulate new multi-disciplinary collaborations.
\vspace{1em}

Finally, we think it is important to emphasize that LLMs are only becoming good at mathematics because they have been \emph{trained} on mathematics. This, once again, proves that the enormous body of research produced by mathematicians is not only a beautiful example of human ingenuity, but is fundamentally useful. In particular, mathematics is -- and always will be -- the language of the universe. With the recent developments in AI, we now have a better way to harness this knowledge, and our objective for the future is best summarized as the desire to use these tools to advance our understanding of mathematics and its implications on the rest of the world.

\section{Acknowledgments}
We thank Alessio Figalli and Joaquim Serra for encouraging us to write this note. We also thank Javier Gomez-Serrano for mentorship and advice as well as several participants at the National Academy of Sciences event \emph{Organizing Mathematical Knowledge in the Age of AI and Formalization}, July 20-21 2026  and the \emph{International Congress of Mathematicians}, July 23-30 2026 for sharing their comments.
Finally, we wish to thank various members of the Banach lattice and phase retrieval communities, including Rima Alaifari, Antonio Avil\'es, Susanna Bertolini, Eugene Bilokopytov, Avik Das, Ingrid Daubechies, Dan Freeman, Enrique García-Sánchez, Josef Greilhuber, David de Hevia, Jorge S.\ Ibáñez-Marcos, Jesús Illescas, Philippe Jaming, Denny Leung, Gonzalo Martínez-Fernández, Timur Oikhberg, Ben Pineau, João Ramos, Christian Rosendal, Alberto Salguero-Alarcón, Alexander Schell, Nazaret Trejo-Arroyo, Vladimir Troitsky and Matthias Wellershoff for the numerous discussions that have influenced the perspectives that we have shared here.
\vspace{1em}

In thanking all of these people, we do not wish to imply that incorporating technology into our mathematical lives has been emphatically desired by all. Change is almost always coupled with uncertainty, and how one adapts to change is a very personal choice. This note only serves to explain what has worked for our community. With this in mind, we are exceptionally grateful to the people mentioned above for the many years of friendship and collaboration, which, in turn, made their suggestions and input over the last few months even more valuable.
\vspace{1em}

The authors declare that they did not use any AI tools when writing this article. The AI models used over the last several months to gather the above experiences were Claude Opus 4.6, 4.7 and 4.8, a small amount of Fable 5, and GPT 5.4 and 5.5. We emphasize that the views and practices described in this paper are not static and will certainly change over time as the models improve.

\bibliographystyle{plain}
\bibliography{sources}

@article{SPRRV,
  title={Stable {P}hase {R}etrieval for {S}pans of {I}ndependent {R}andom {V}ariables},
  author={Abdalla, Pedro and de Dios Pont, Jaume and Ramos, Jo\~ao P.~G. and Taylor, Mitchell~A.},
  journal={arXiv preprint arXiv:2607.06693},
  year={2026}
}

@article{Lukas1,
  title={On the existence problem of regular {G}abor frames},
  author={Liehr, Lukas and de Dios Pont, Jaume and Taylor, Mitchell~A.},
  journal={arXiv preprint arXiv:2606.26052},
  year={2026}
}

@article{Lukas2,
  title={Cantor measures with odd base do not admit {F}ourier frames},
  author={Liehr, Lukas and de Dios Pont, Jaume and Taylor, Mitchell~A.},
  journal={arXiv preprint arXiv:2607.08656},
  year={2026}
}

@article{bertolini20262,
  title={{$L^2$}-{S}tability for {STFT} phase retrieval},
  author={Bertolini, Susanna and Pont, Jaume de Dios and Pineau, Ben and Taylor, Mitchell A. and Ramos, Jo{\~a}o P.~G.},
  journal={arXiv preprint arXiv:2605.20527},
  year={2026}
}

@article{jiang2026,
      title={From Solvers to Research: Large Language Model-Driven Formal Mathematics at the Research Frontier}, 
      author={Eric Jiang and Xiao Liang and Yikai Zhang and Yingjia Wan and Mengting Li and Haikang Deng and Alexander K. Taylor and Justin Baker and Rushil Raghavan and Junyi Zhang and Ying Nian Wu and Andrea L. Bertozzi and Kai-Wei Chang and Raghu Meka and Matthew Sottile and Nanyun Peng and Amit Sahai and Terence Tao and Wei Wang},
      year={2026},
      journal={arXiv preprint arXiv:2607.07779}
}

@inproceedings{mathlib2020,
  author    = {{The mathlib {C}ommunity}},
  title     = {The {L}ean {M}athematical {L}ibrary},
  booktitle = {Proceedings of the 9th {ACM} {SIGPLAN} International Conference
               on Certified Programs and Proofs},
  series    = {CPP 2020},
  publisher = {ACM},
  address   = {New Orleans, LA, USA},
  year      = {2020},
  month     = jan,
  doi       = {10.1145/3372885.3373824},
  url       = {https://doi.org/10.1145/3372885.3373824}
}

@inproceedings{lean4,
  title = {The {L}ean 4 {T}heorem {P}rover and {P}rogramming {L}anguage},
  author = {de Moura, Leonardo and Ullrich, Sebastian},
  year = {2021},
  isbn = {978-3-030-79875-8},
  publisher = {Springer-Verlag},
  address = {Berlin, Heidelberg},
  url = {https://doi.org/10.1007/978-3-030-79876-5_37},
  doi = {10.1007/978-3-030-79876-5_37},
  booktitle = {Automated Deduction – CADE 28: 28th International Conference on Automated Deduction, Virtual Event, July 12–15, 2021, Proceedings},
  pages = {625–635},
  numpages = {11}
}

@article{GT,
 author = {Garc{\'{\i}}a-S{\'a}nchez, E. and Tradacete, P.},
 title = {Free dual spaces and free {Banach} lattices},
 fjournal = {Journal of Mathematical Analysis and Applications},
 journal = {J. Math. Anal. Appl.},
 issn = {0022-247X},
 volume = {532},
 number = {2},
 pages = {22},
 note = {Id/No 127931},
 year = {2024},
 language = {English},
 doi = {10.1016/j.jmaa.2023.127931},
 zbMATH = {7787735},
 Zbl = {1534.46014}
}

@article{AT,
 author = {Avil{\'e}s, Antonio and Tradacete, Pedro},
 title = {Amalgamation and injectivity in {Banach} lattices},
 fjournal = {IMRN. International Mathematics Research Notices},
 journal = {Int. Math. Res. Not.},
 issn = {1073-7928},
 volume = {2023},
 number = {2},
 pages = {956--997},
 year = {2023},
 language = {English},
 doi = {10.1093/imrn/rnab285},
 zbMATH = {7652752},
 Zbl = {1518.46010}
}

@misc{MFT,
 author = {Gonzalo Mart{\'{\i}}nez-Fern{\'a}ndez and Pedro Tradacete},
 title = {Free {Products} of {Banach} {Lattices}},
 year = {2026},
 howpublished = {\textit{arXiv preprint {arXiv}:2605.28988}},
 url = {https://arxiv.org/abs/2605.28988},
 arXiv = {arXiv:2605.28988}
}

@misc{tao2026,
  author       = {Terence Tao},
  title        = {The Integrated Explicit Analytic Number Theory Network : a progress report},
  howpublished = {Conference talk at Techniques and Tools for the Formalization of Analysis, Institute for Computational and Experimental Research in Mathematics (ICERM)},
  year         = {2026},
  url          = {https://icerm.brown.edu/video_archive/4649},
  urldate      = {2026-07-03}
}

@article{jerison2000hot,
  title={The “hot spots” conjecture for domains with two axes of symmetry},
  author={Jerison, David and Nadirashvili, Nikolai},
  journal={Journal of the American Mathematical Society},
  volume={13},
  number={4},
  pages={741--772},
  year={2000}
}

@article{ilin2026sorries,
  title={{Sorries Are Not the Hard Part: An Expert-Review Case Study of a Semi-Autonomous Formalization}},
  author={Ilin, Vasily and Nugent, Brian},
  journal={arXiv preprint arXiv:2606.13925},
  year={2026}
}

@article{steinerberger2023upper,
  title={An upper bound on the hot spots constant},
  author={Steinerberger, Stefan},
  journal={Revista Mathematica Iberoamericana},
  volume={39},
  number={4},
  year={2023}
}

@article{mariano2023improved,
  title={Improved upper bounds for the {H}ot {S}pots constant of {L}ipschitz domains},
  author={Mariano, Phanuel and Panzo, Hugo and Wang, Jing},
  journal={Potential Analysis},
  volume={59},
  number={2},
  pages={771--787},
  year={2023},
  publisher={Springer}
}

@article{pont2025sharp,
  title={Sharp bounds on the failure of the hot spots conjecture},
  author={Pont, Jaume de Dios and Hsu, Alexander W and Taylor, Mitchell A},
  journal={arXiv preprint arXiv:2508.16321},
  year={2025}
}

@article{rammal2026formalizing,
  title={Formalizing mathematics at scale},
  author={Rammal, Ahmad and Patel, Niket and Gloeckle, Fabian and Hayat, Amaury and Kempe, Julia and Munos, Remi and Arnal, Charles and Cabannes, Vivien},
  journal={arXiv preprint arXiv:2605.29955},
  year={2026}
}

@incollection{thurston2006proof,
  title={On proof and progress in mathematics},
  author={Thurston, William P},
  booktitle={18 Unconventional essays on the nature of mathematics},
  pages={37--55},
  year={2006},
  publisher={Springer}
}

@article{avigad2026mathematicians,
  title={{Mathematicians in the Age of AI}},
  author={Avigad, Jeremy},
  journal={arXiv preprint arXiv:2603.03684},
  year={2026}
}

@article{klowden2026mathematical,
  title={{Mathematical methods and human thought in the age of AI}},
  author={Klowden, Tanya and Tao, Terence},
  journal={arXiv preprint arXiv:2603.26524},
  year={2026}
}

@misc{BSTT,
 author = {E.~Bilokopytov and A.~Salguero-Alarc\'on and P.~Tradacete and N.~Trejo-Arroyo},
 title = {Splitting and extension properties in {B}anach lattices},
 howpublished = {Work in progress}
}

@article{BanLib,
 author = {David~Muñoz-Lahoz},
 title = {The {B}anach lattice {L}ean library},
 year = {2026}
}

@article{davies2021advancing,
  title={Advancing mathematics by guiding human intuition with AI},
  author={Davies, Alex and Veli{\v{c}}kovi{\'c}, Petar and Buesing, Lars and Blackwell, Sam and Zheng, Daniel and Toma{\v{s}}ev, Nenad and Tanburn, Richard and Battaglia, Peter and Blundell, Charles and Juh{\'a}sz, Andr{\'a}s and others},
  journal={Nature},
  volume={600},
  number={7887},
  pages={70--74},
  year={2021},
  publisher={Nature Publishing Group UK London}
}

@inproceedings{buzzard2020formalising,
  title={Formalising perfectoid spaces},
  author={Buzzard, Kevin and Commelin, Johan and Massot, Patrick},
  booktitle={Proceedings of the 9th ACM SIGPLAN International Conference on Certified Programs and Proofs},
  pages={299--312},
  year={2020}
}

@article{tao2025machine,
  title={Machine-assisted proof},
  author={Tao, Terence},
  journal={Notices of the American Mathematical Society},
  volume={72},
  number={1},
  pages={6--13},
  year={2025},
  publisher={American Mathematical Society}
}

@misc{leantutorial,
 author = {David~Muñoz-Lahoz},
 title = {Autoformalize your math},
 note   = {\url{https://lean.functionalanalysismadrid.com/}},
}

@article{scholze2022liquid,
  title={Liquid tensor experiment},
  author={Scholze, Peter},
  journal={Experimental Mathematics},
  volume={31},
  number={2},
  pages={349--354},
  year={2022},
  publisher={Taylor \& Francis}
}

@article{OTTT,
 author = {Oikhberg, Timur and Taylor, Mitchell A. and Tradacete, Pedro and Troitsky, Vladimir G.},
 title = {Free {Banach} lattices},
 fjournal = {Journal of the European Mathematical Society (JEMS)},
 journal = {J. Eur. Math. Soc. (JEMS)},
 issn = {1435-9855},
 volume = {28},
 number = {10},
 pages = {4387--4514},
 year = {2026},
 language = {English},
 doi = {10.4171/jems/1552},
 zbMATH = {8221493}
}

\end{document}